\documentclass[a4paper,11pt]{article}
\usepackage[left=25mm,right=20mm,top=20mm,bottom=20mm]{geometry}
\usepackage{amsmath}
\usepackage{amsfonts}
\usepackage{amssymb}
\usepackage{amsbsy}
\usepackage{graphics}
\usepackage{setspace}
\usepackage{epsfig}
\usepackage{cite}
\usepackage{graphicx}
\usepackage{latexsym}
\newtheorem{theorem}{Theorem}[section]
\newtheorem{lemma}[theorem]{Lemma}

\newtheorem{definition}[theorem]{Definition}
\newtheorem{example}[theorem]{Example}

\def\bproof{\noindent\textbf{Proof:\ }}
\def\eproof{$\Box$\\}
\numberwithin{equation}{section}
\begin{document}
\title{Guaranteed upper and lower bounds on the uniform load of contact problems in elasticity}
\author{L. Angela Mihai\footnote{School of Mathematics, Cardiff University, Senghennydd Road, Cardiff, CF24 4AG, UK, Email: MihaiLA@cardiff.ac.uk} \qquad Alain Goriely\footnote{Mathematical Institute, University of Oxford, Woodstock Road, Oxford, OX2 6GG, UK, Email: goriely@maths.ox.ac.uk}}\date{}
\maketitle

\begin{abstract}
Two mathematical models are developed within the theoretical framework of large strain elasticity for the determination of upper and lower bounds on the total strain energy of a finitely deformed hyperelastic body in unilateral contact with a rigid surface or with an elastic substrate. The model problems take the form of two continuous optimization problems with inequality constraints, the solutions of which are used to provide an enclosure on the uniform external load acting on the body's surface away from the contact zone. \\

\end{abstract}

\section{Introduction}
In materials analysis and design, contact problems in elasticity are central to the modelling and investigation of many structural systems. In particular, the analysis of soft tissue biomechanics and the design of bioinspired synthetic structures involve non-linear hyperelastic models for which the mathematical and numerical treatment poses many physical, theoretical, and computational challenges \cite{Auricchio:2013:AVLRTW,Hueber:2012:HW,Mihai:2011:MG,Mihai:2015:MAG,Mihai:2015:MCJG,Nedoma:2003:NKFKS,Sitzmann:2014:SWW}.

Hyperelastic materials are the class of material models described by a strain energy density function \cite{Green:1970:GA,Green:1968:GZ,Ogden:1997,TruesdellNoll:2004}. For these materials, boundary value problems are often equivalent to variational problems, which provide powerful methods for obtaining approximate solutions. They can also be used to generate finite element methods for which the numerical analysis stands on the shoulder of the mathematical analysis for the elastic model \cite{Carstensen:2004:CD,LeTallec94,Mihai:2013:MG,Oden:1972}.

The complexity of contact problems modelling structural systems is generally associated with the detection of contacts and openings and the resolution of non-linear equations for contact. For example, in systems formed from linear elastic bodies (\emph{e.g.} ceramics, metals), surface roughness can impede active contact when a surface is pressed against another, and contact forces cannot be transmitted where surface separation occurs \cite{Hlavacek:1988:HHNL,Johnson:1985,Kikuchi:1988:KO}. By contrast, non-linear elastic bodies (\emph{e.g.} rubber, soft tissue) are more pliable and thus capable of attaining more active support through which contact forces can be transmitted effectively \cite{Ball:2002,Ciarlet:1985:CN,Ciarlet:1987:CN,Johnson:1971:JKR}. These problems can be formulated and solved in the framework of variational inequalities \cite{Lions:1967:LS}, which originated in the paper by Fichera \cite{Fichera:1964} on the existence and uniqueness of solution to the celebrated \emph{Signorini problem} with \emph{ambiguous boundary conditions} \cite{Signorini:1933}. While Signorini problem consisted in finding the equilibrium state of an elastic body resting on a rigid frictionless surface, the influence of variational inequalities went beyond the fields of mathematical analysis, optimization, and mechanics \cite{Duvaut:1976:DL,Glowinski:1981:GLT,Wriggers:2006}, leading also to new models in game theory, economics, and finance.

In large strain elasticity, the study of finite dimensional structural models is motivated by the increasing need for effective computational techniques required in practical applications, but little attention has been paid to date to the underpinning infinite dimensional problems which are fundamental in the development of new methods. For example, while the reciprocal theorems of stationary potential and complementary energy are well known in the infinitesimal theory of elasticity and in the theory of structures \cite{Christiansen:1996,Mihai:2009:MA}, for large strain deformations, the stress-strain relation is non-linear, and  its inverse is not uniquely defined in general, either locally or globally. For systems with a finite number of degrees of freedom, extensions of the reciprocal principles have been proposed in \cite{Langhaar:1953,Libove:1964,Pipes:1962}. For an elastic continuum, assuming that the stress-strain relation is invertible, a complementary energy principle in terms of stresses was formulated via a Legendre transformation on the strain energy function in \cite{Levinson:1965}. In order to avoid the difficulty of inverting the constitutive relation, in \cite{Lee:1980:LS,Shield:1980:SL}, trial functions for the deformation gradient instead of the stress tensor are used to obtain a variational principle of the complementary energy type. Then, under appropriate conditions, the stationary potential energy and complementary energy principles become extremum principles which can be used to provide lower and upper bounds on physical quantities of interest.

In the present study, the variational approach of \cite{Lee:1980:LS} is extended first to the problem of a non-linear hyperelastic body in unilateral contact with a rigid surface, then to the case of two hyperelastic bodies in mutual non-penetrative contact. For each of these problems, the corresponding variational models take the form of two continuous optimization problems with inequality constraints, the solutions of which can be used to provide upper and lower bounds for the uniform external loading acting away from the contact zone. To illustrate the theory, analytical upper and lower bounds for the external load of a system formed from two elastic bodies in mutual unilateral contact and subject to large compression or bending, which can be maintained in every homogeneous isotropic incompressible hyperelastic material in the absence of body forces, are obtained.

In practice, contact problems with large stresses and strains at the adjoining material surfaces arise, for example, when creases are formed and self-contact takes place in soft materials and structures \cite{Zalachas:2013:ZCSL}, or in biological systems where the attachment between cells are sufficiently weak so that cells separate, and failure through the appearance of gaps between cells occurs \cite{Lewis:2008:LYMC}. Although cell debonding is a spontaneous mechanism for crack initiation in many natural structures, it has been less investigated to date \cite{Bassel:2014:etal,Bruce:2003,Gao:1990:GPR}.

The remainder of this paper is organised as follows: in Section~\ref{CC:sec:problem:one}, the equilibrium problem is formulated for a non-linear elastic body in non-penetrative contact with itself or with a rigid obstacle, in the absence of friction forces; in Sections~\ref{CC:sec:bounds} and \ref{CC:sec:bounds:con}, variational problems are derived and analysed for the upper and lower bounds on the total strain energy of a finitely deformed elastic body made from a compressible or an incompressible material, respectively; in Section~\ref{CC:sec:problem:two}, the variational approach is extended to the case of two elastic bodies in mutual non-penetrative contact, and examples where the elastic bodies are made from a neo-Hookean material and subject to large compression or bending are presented.

\section{Elastostatic Equilibrium with Unilateral Contact}\label{CC:sec:problem:one}
\begin{figure}[htbp]
\begin{center}
\scalebox{0.5}{\includegraphics{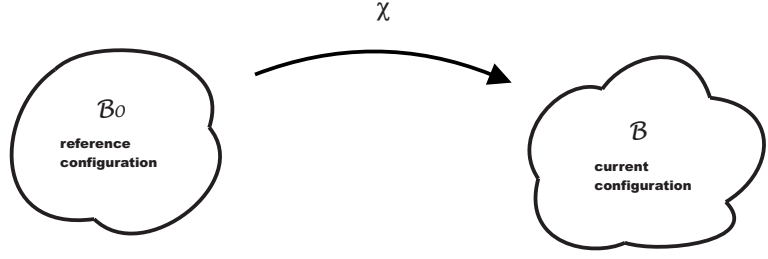}}
\caption{Schematic representation of elastic deformation.}\label{fig:deformation}
\end{center}
\end{figure}

A continuous material body occupies a compact domain $\bar\Omega$ of the three-dimensional Euclidean space $\mathbb{R}^{3}$, such that the interior of the body is an open, bounded, connected set $\Omega\subset\mathbb{R}^{3}$, and its boundary $\Gamma=\partial\Omega=\bar\Omega\setminus\Omega$ is Lipschitz continuous (in particular, we assume that a unit normal vector $\textbf{n}$ exists almost everywhere on $\Gamma$). The body is subject to a finite elastic deformation defined by the one-to-one, orientation preserving transformation:
\[
\boldsymbol{\chi}:\Omega\to\mathbb{R}^{3},
\]
such that  $J=\det\left(\mathrm{Grad}\boldsymbol{\chi}\right)>0$ on $\Omega$  and $\boldsymbol{\chi}$ is injective on $\Omega$ (see Figure~\ref{fig:deformation}). The injectivity condition on $\Omega$ guarantees that interpenetration of the matter is avoided. However, since self-contact is permitted, this transformation does not need to be injective on $\bar\Omega$.

Let the spatial point $\textbf{x}=\boldsymbol{\chi}(\textbf{X})$ correspond to the place occupied by the particle $\textbf{X}$ in the deformation $\boldsymbol{\chi}$. For the deformed body, the equilibrium state in the absence of a body load is described in terms of the Cauchy stress by the Eulerian field equation:
\[
\mathrm{div}\ \boldsymbol{\sigma}(\textbf{x})=\textbf{0}.
\]
The above governing equation is completed by a constitutive law for $\boldsymbol{\sigma}$, depending on material properties, and supplemented by boundary conditions.

Since the domain occupied by the body after deformation is usually unknown, we rewrite the above equilibrium problem as an equivalent problem in the reference configuration where the independent variables are $\textbf{X}\in\Omega$. The corresponding Lagrangian equation of non-linear elastostatics is:
\begin{eqnarray}
\mathrm{Div}\ \textbf{P}(\textbf{X})=\textbf{0} &\mbox{in}& \Omega,\label{CC:eq:balance}
\end{eqnarray}
where $\textbf{P}$ is the first Piola-Kirchhoff stress tensor.

For a homogeneous compressible hyperelastic material described by the strain energy function $\mathcal{W}(\textbf{F})$, the first Piola-Kirchhoff stress tensor is:
\begin{eqnarray}
\textbf{P}=\frac{\partial\mathcal{W}}{\partial\textbf{F}}&\mbox{in}& \Omega,\label{CC:eq:Pstress}
\end{eqnarray}
where $\textbf{F}=\mathrm{Grad}\ \boldsymbol{\chi}$ is the deformation gradient.

Then the corresponding Cauchy stress tensor can be expressed as follows
$\boldsymbol{\sigma}=J^{-1}\textbf{P}\textbf{F}^{T}$, and $\textbf{P}=\boldsymbol{\sigma}\mathrm{cof}\textbf{F}$.

\begin{figure}[htbp]
\begin{center}
\scalebox{0.5}{\includegraphics{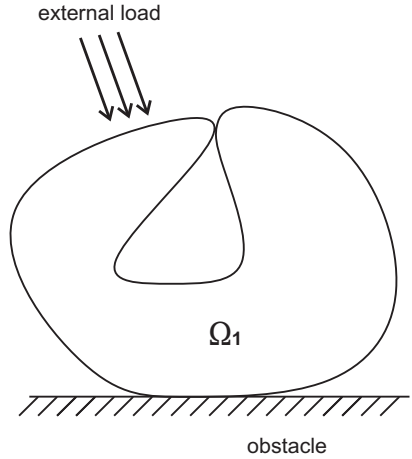}}
\caption{Schematic representation of an elastic body in non-penetrative contact with itself and with a rigid obstacle and subject to external load.}\label{fig:contact}
\end{center}
\end{figure}

If the body is in non-penetrative contact with itself or with a rigid obstacle and subject to external dead loading conditions (see Figure~\ref{fig:contact}), then the general boundary value problem is to find the displacement field $\textbf{u}(\textbf{X})=\textbf{x}-\textbf{X}$, for all $\textbf{X}\in\Omega$, such that the equilibrium equation (\ref{CC:eq:balance}) is satisfied subject to the following conditions on the relatively disjoint, open subsets of the boundary $\{\Gamma_{D},\Gamma_{N},\Gamma_{C}\}\subset\partial\Omega$, such that $\partial \Omega\setminus\left(\Gamma_{D}\cup\Gamma_{N}\cup\Gamma_{C}\right)$ has zero area \cite[Chapter VI]{LeTallec94}:
\begin{itemize}
\item On $\Gamma_{D}$, the Dirichlet (displacement) conditions :
\begin{eqnarray}
\textbf{u}(\textbf{X})=\textbf{u}_{D},\label{CC:eq:Dbc}
\end{eqnarray}
\item On $\Gamma_{N}$, the Neumann (traction) conditions :
\begin{eqnarray}
\textbf{P}(\textbf{X})\textbf{N}=\textbf{g}_{N}(\textbf{X}),\label{CC:eq:Nbc}
\end{eqnarray}
where $\textbf{N}$ is the outward unit normal vector to $\Gamma_{N}$, and $\textbf{g}_{N}dA=\boldsymbol{\tau}da$, where $\boldsymbol{\tau}=\boldsymbol{\sigma}\textbf{n}$ is the surface traction measured per unit area of the deformed state.

\item On $\Gamma_{C}$, the conditions for unilateral (non-penetrative) contact either with itself or with a rigid obstacle, in the absence of friction forces:
\begin{eqnarray}
\eta(\textbf{X}+\textbf{u}(\textbf{X}))&\leq& d,\label{CC:eq:ndcc}\\
\textbf{P}(\textbf{X})\textbf{N}\cdot\textbf{N}=\textbf{g}_{C}(\textbf{X})\cdot\textbf{N}&\leq& g,\label{CC:eq:nscc}\\
\left(\eta(\textbf{X}+\textbf{u}(\textbf{X}))-d\right)\left(\textbf{g}_{C}(\textbf{X})\cdot\textbf{N}-g\right)&=&0,\label{CC:eq:cpcc}\\
\textbf{g}_{C}(\textbf{X})=\textbf{P}(\textbf{X})\textbf{N}=-\textbf{P}(\textbf{X}')\textbf{N}'=-\textbf{g}_{C}(\textbf{X}') &\mbox{if}& \boldsymbol{\chi}(\textbf{X})=\boldsymbol{\chi}(\textbf{X}').\label{CC:eq:sccc}
\end{eqnarray}
For the contact with a rigid obstacle, $\eta:\mathbb{R}^3\to\mathbb{R}$ is the function describing the relative position of the body to the surface of the obstacle and $\textbf{N}$ is the outward unit normal vector to this surface (oriented towards the obstacle), $d\geq0$ is the relative distance which cannot be exceeded between potential contact points, and $g\geq0$ is the cohesion parameter. For self-contact, $\eta=[\textbf{u}(\textbf{X})]\cdot\textbf{N}$ where $[\textbf{u}(\textbf{X})]=\textbf{u}(\textbf{X})-\textbf{u}(\textbf{X}')$ denotes the jump in the displacement at the points $\textbf{X}\neq\textbf{X}'$ where self-contact may occur and $\textbf{N}$ and $\textbf{N}'$ represent the outward unit normal vectors to the boundary at those points, respectively. Then (\ref{CC:eq:ndcc}) sets the permitted relative distance between potential contact points; (\ref{CC:eq:nscc}) gives the allowed normal force acting at a contact point; (\ref{CC:eq:cpcc}) is the complementarity condition that the maximum relative distance or the maximum normal force is attained at each point; and (\ref{CC:eq:sccc}) states that, at the points of self-contact, the action and reaction principle holds. 
\end{itemize}

For the study of the above contact boundary value problem (CBVP), the following result guarantees injectivity of the deformation mapping, so that self-penetration of the elastic body cannot occur \cite{Ball:2002}.

\begin{theorem}
If $\boldsymbol{\chi}\in C^{1}(\bar{\Omega},\mathbb{R}^{3})$, such that:
\begin{eqnarray}
\mathrm{J}=\det(\mathrm{Grad}\ \boldsymbol{\chi})>0 &\mbox{on}& \Omega,\label{CC:eq:inj1}
\end{eqnarray}
satisfies the following condition:
\begin{eqnarray}
\int_{\Omega}\det(\mathrm{Grad}\ \boldsymbol{\chi})d\mathrm{\textbf{X}}&\leq&\mathrm{vol}\left(\boldsymbol{\chi}(\Omega)\right),\label{CC:eq:inj2}
\end{eqnarray}
then $\boldsymbol{\chi}$ is injective on $\Omega$.
\end{theorem}

\bproof (See \emph{e.g.} \cite[p.~571]{LeTallec94})
Assuming that $\boldsymbol{\chi}$ is not injective on $\Omega$, it follows that there exists an open ball $\mathcal{B}\subset\boldsymbol{\chi}(\Omega)$ whose elements are the image of at least two points of $\Omega$, \emph{i.e.}:
\[
\mathrm{card}\left(\boldsymbol{\chi}(\Omega)\right)\geq2,\qquad \forall\textbf{X}\in\mathcal{B}.
\]
Therefore:
\[
\mathrm{vol}\left(\boldsymbol{\chi}(\Omega)\right)=\int_{\boldsymbol{\chi}(\Omega)}d\mathrm{\textbf{X}}<\int_{\boldsymbol{\chi}(\Omega)}\mathrm{card}\left(\boldsymbol{\chi}(\Omega)\right)d\mathrm{\textbf{X}}=\int_{\Omega}\det(\mathrm{Grad}\ \boldsymbol{\chi})d\mathrm{\textbf{X}},
\]
which contradicts (\ref{CC:eq:inj2}).
\eproof

\section{Unconstrained Materials}\label{CC:sec:bounds}

We now cast the CBVP (\ref{CC:eq:balance})-(\ref{CC:eq:sccc}) in variational forms similar to those described in \cite{Lee:1980:LS}. However, in our case, we also take into account the kinematic constraints on the displacement field and the associated reaction force (\ref{CC:eq:ndcc})-(\ref{CC:eq:nscc}) at the potential contact zone $\Gamma_{C}$. Throughout the analysis, we assume that the conditions (\ref{CC:eq:inj1})-(\ref{CC:eq:inj2}) are satisfied. For simplicity of notation, we also set $d=0$ and $g=0$ in (\ref{CC:eq:ndcc})-(\ref{CC:eq:cpcc}), which corresponds to the physical case of non-penetrative cohesionless contact. The case when these parameters are non-zero can be treated by analogy. For the given problem, the existence of a solution implies the mathematical correctness of the problem, while the (local) uniqueness of the solution indicates that the problem can be used to model of a physical phenomenon.

\begin{definition} (kinematically admissible fields) The vector field $\textbf{u}$ is a kinematically admissible displacement if and only if it satisfies simultaneously the Dirichlet boundary condition (\ref{CC:eq:Dbc}) and the non-penetrative contact condition (\ref{CC:eq:ndcc}).
\end{definition}

We denote the closed convex set of kinematically admissible displacement fields by:
\[
\mathcal{K}=\left\{\textbf{u}'\in W^{1,s}(\Omega;\mathbb{R}^3)\ |\ \textbf{u}'(\textbf{X})=\textbf{u}_{D} \ \mbox{on} \ \Gamma_{D},\ \eta(\textbf{X}+\textbf{u}'\left(\textbf{X})\right)\leq 0\ \mbox{on} \ \Gamma_{C}\right\},
\]
for some $s>3/2$ (see \emph{e.g.} \cite{LeTallec94}).

\begin{definition} (statically admissible fields)
A tensor field $\textbf{F}$ is a statically admissible gradient if and only if $\textbf{P}$ given by (\ref{CC:eq:Pstress}) satisfies simultaneously the equilibrium condition (\ref{CC:eq:balance}), the boundary condition of dead loading (\ref{CC:eq:Nbc}), and the contact conditions (\ref{CC:eq:nscc}) and (\ref{CC:eq:sccc}).
\end{definition}

The closed convex set of statically admissible gradient fields is denoted by:
\[
\begin{split}
\mathcal{S}=\{\textbf{F}'\in L^{q}(\Omega;\mathbb{R})\ |\ &\boldsymbol{\mathrm{Div}}\frac{\partial\mathcal{W}}{\partial\textbf{F}'}=\textbf{0}\ \mbox{in}\ \Omega, \\
&\frac{\partial\mathcal{W}}{\partial\textbf{F}'}\textbf{N}=\textbf{g}_{N}\ \mbox{on} \ \Gamma_{N}, \\
&\frac{\partial\mathcal{W}}{\partial\textbf{F}'}\textbf{N}\cdot\textbf{N}=\textbf{g}_{C}\cdot\textbf{N}\leq0 \ \mbox{on} \ \Gamma_{C},\\
&\textbf{g}_{C}(\textbf{X})=-\textbf{g}_{C}(\textbf{X}') \ \mbox{if}\ \boldsymbol{\chi}(\textbf{X})=\boldsymbol{\chi}(\textbf{X}') \ \mbox{on} \ \Gamma_{C}\},
\end{split}
\]
for some $q>1$ (see \cite{LeTallec94}).

Let $\textbf{u}',\textbf{u}''\in\mathcal{K}$ and $\textbf{F}''\in\mathcal{S}$, such that $\textbf{F}''=\textbf{I}+\mathrm{Grad}\ \textbf{u}''$. By the divergence theorem, we obtain:
\begin{equation}\label{CC:eq:mixed}
\begin{split}
\int_{\Gamma_{D}}\frac{\partial\mathcal{W}}{\partial\textbf{F}''}\textbf{N}\cdot\textbf{u}_{D}dA+\int_{\Gamma_{N}}\textbf{g}_{N}\cdot\textbf{u}'dA+\int_{\Gamma_{C}}\frac{\partial\mathcal{W}}{\partial\textbf{F}''}\textbf{N}\cdot\textbf{u}'dA
&=\int_{\Omega}\mathrm{Div}\left(\frac{\partial\mathcal{W}}{\partial\textbf{F}''}\textbf{u}'\right)dV\\
&=\int_{\Omega}\frac{\partial\mathcal{W}}{\partial\textbf{F}''}:\textbf{F}'dV,
\end{split}
\end{equation}
where $\textbf{F}'=\textbf{I}+\mathrm{Grad}\ \textbf{u}'$. From this identity we derive two important inequalities as follows.

If, in (\ref{CC:eq:mixed}), $\textbf{u}'=\textbf{u}''=\textbf{u}\in\mathcal{K}$ and $\textbf{F}''=\textbf{F}\in\mathcal{S}$ are the displacement and gradient fields, respectively, satisfying the CBVP, then:
\begin{eqnarray}
\int_{\Gamma_{D}}\frac{\partial\mathcal{W}}{\partial\textbf{F}}\textbf{N}\cdot\textbf{u}_{D}dA+\int_{\Gamma_{N}}\textbf{g}_{N}\cdot\textbf{u}dA
&=&\int_{\Omega}\frac{\partial\mathcal{W}}{\partial\textbf{F}}:\textbf{F}dV.\label{CC:eq:umixed}
\end{eqnarray}

\begin{enumerate}
\item First, in (\ref{CC:eq:mixed}), we take $\textbf{u}'\in\mathcal{K}$ and set $\textbf{u}''=\textbf{u}\in\mathcal{K}$ and $\textbf{F}''=\textbf{F}\in\mathcal{S}$ to satisfy the CBVP. Subtracting (\ref{CC:eq:umixed}) from (\ref{CC:eq:mixed}) implies:
\begin{eqnarray*}
\int_{\Gamma_{N}}\textbf{g}_{N}\cdot(\textbf{u}'-\textbf{u})dA+\int_{\Gamma_{C}}\frac{\partial\mathcal{W}}{\partial\textbf{F}}\textbf{N}\cdot\textbf{u}'dA
&=&\int_{\Omega}\frac{\partial\mathcal{W}}{\partial\textbf{F}}:(\textbf{F}'-\textbf{F})dV,
\end{eqnarray*}
where the second integral is non-negative.

The primal (displacement) variational problem is to find $\textbf{u}\in\mathcal{K}$ satisfying:
\begin{equation}\label{CC:eq:balance:vdp}
\int_{\Omega}\frac{\partial\mathcal{W}}{\partial\textbf{F}}\left(\textbf{X},\textbf{F}(\textbf{u})\right): \mathrm{Grad}(\textbf{u}'-\textbf{u})dV
\geq\int_{\Gamma_{N}}\textbf{g}_{N}\cdot(\textbf{u}'-\textbf{u})dA,
\end{equation}
for all admissible fields $\textbf{u}'\in\mathcal{K}$.

\item Next, in (\ref{CC:eq:mixed}), we take $\textbf{F}''\in\mathcal{S}$ and set $\textbf{u}'=\textbf{u}\in\mathcal{K}$ and $\textbf{F}'=\textbf{F}\in\mathcal{S}$ to satisfy the CBVP. Subtracting (\ref{CC:eq:umixed}) from (\ref{CC:eq:mixed}) yields:
\begin{eqnarray*}
\int_{\Gamma_{D}}\left(\frac{\partial\mathcal{W}}{\partial\textbf{F}''}-\frac{\partial\mathcal{W}}{\partial\textbf{F}}\right)\textbf{N}\cdot\textbf{u}_{D}dA
+\int_{\Gamma_{C}}\frac{\partial\mathcal{W}}{\partial\textbf{F}''}\textbf{N}\cdot\textbf{u}dA
&=&\int_{\Omega}\left(\frac{\partial\mathcal{W}}{\partial\textbf{F}''}-\frac{\partial\mathcal{W}}{\partial\textbf{F}}\right):\textbf{F}dV,
\end{eqnarray*}
where the second integral is again non-negative.

The complementary-type variational problem is to find $\textbf{F}\in\mathcal{S}$ such that:
\begin{equation}\label{CC:eq:balance:vgp}
\int_{\Omega}\left(\frac{\partial\mathcal{W}}{\partial\textbf{F}'}-\frac{\partial\mathcal{W}}{\partial\textbf{F}}\right):\textbf{F}dV
\geq\int_{\Gamma_{D}}\left(\frac{\partial\mathcal{W}}{\partial\textbf{F}'}-\frac{\partial\mathcal{W}}{\partial\textbf{F}}\right)\textbf{N}\cdot\textbf{u}_{D}dA,
\end{equation}
for all admissible fields $\textbf{F}'\in\mathcal{S}$.
\end{enumerate}

Assuming that the contact problem CBVP has a statically admissible solution $\textbf{F}\in\mathcal{S}$ and a corresponding kinematically admissible displacement $\textbf{u}\in\mathcal{K}$, our objective is to determine upper and lower bounds on the total strain energy of the deformed body.

\subsection{Upper and Lower Bounds on the Strain Energy}
The potential energy of a kinematically admissible vector field $\textbf{u}'\in\mathcal{K}$ is defined by:
\begin{eqnarray}
E_{p}(\textbf{u}')&=&\int_{\Omega}\mathcal{W}\left(\textbf{F}'(\textbf{u}')\right)dV-\int_{\Gamma_{N}}\textbf{g}_{N}\cdot\textbf{u}'dA.\label{CC:eq:penergy}
\end{eqnarray}
Note that the contact constraint appears in the set of admissible fields $\mathcal{K}$, but not in the expression of the energy functional $E_{p}$.

The \emph{principle of stationary potential energy} states that the solution of the CBVP is characterised by the variational formulation:
\[
\mbox{find} \quad \textbf{u}\in \mathcal{K} \quad \mbox{such that}\quad
E_{p}(\textbf{u})=\inf_{\textbf{u}'\in\mathcal{K}}E_{p}(\textbf{u}').
\]

\begin{lemma}
The energy functional $E_{p}$ defined by (\ref{CC:eq:penergy}) has a local minimum at
the solution $\textbf{u}\in\mathcal{K}$ of the primal variational problem (\ref{CC:eq:balance:vdp}),
\emph{i.e.}
\begin{eqnarray}
E_{p}(\textbf{u})&\leq&E_{p}(\textbf{u}'),\label{CC:eq:ubound}
\end{eqnarray}
for all $\textbf{u}'\in\mathcal{K}$.
\end{lemma}

\bproof Let $\textbf{u}'=\textbf{u}+\boldsymbol{\epsilon}\in\mathcal{K}$, where $\textbf{u}\in\mathcal{K}$ is the solution to (\ref{CC:eq:balance:vdp}) and $\epsilon_{i}\ll 1$, $i=1,2,3$, such that:
\begin{eqnarray}
\boldsymbol{\epsilon}=0\ \mbox{on}\  \Gamma_{D}\cup\Gamma_{C}.\label{CC:eq:epsilon:1}
\end{eqnarray}
Then, the value of the potential energy at $\textbf{u}'$ can be approximated to the second order in $\boldsymbol{\epsilon}$ as follows:
\[
\begin{split}
E_{p}(\textbf{u}')
=&E_{p}(\textbf{u})
+\int_{\Omega}\frac{\partial\mathcal{W}}{\partial F_{ij}}\left(\textbf{X},\textbf{F}(\textbf{u})\right)\frac{\partial\epsilon_{i}}{\partial X_{j}}dV-\int_{\Gamma_{N}}\textbf{g}_{N}\cdot\boldsymbol{\epsilon}dA\\
&+\frac{1}{2}\int_{\Omega}\frac{\partial^{2}\mathcal{W}}{\partial F_{kl}\partial F_{ij}}\left(\textbf{X},\textbf{F}(\textbf{u})\right)\frac{\partial\epsilon_{i}}{\partial X_{j}}\frac{\partial\epsilon_{k}}{\partial X_{l}}dV.
\end{split}
\]
By (\ref{CC:eq:balance:vdp}):
\begin{eqnarray*}
\int_{\Omega}\frac{\partial\mathcal{W}}{\partial F_{ij}}\left(\textbf{X},\textbf{F}(\textbf{u})\right)\frac{\partial\epsilon_{i}}{\partial X_{j}}dV&\geq&\int_{\Gamma_{N}}\textbf{g}_{N}\cdot\boldsymbol{\epsilon}dA.
\end{eqnarray*}
Hence:
\begin{eqnarray*}
E_{p}(\textbf{u}')&\geq&E_{p}(\textbf{u})+
\frac{1}{2}\int_{\Omega}\frac{\partial^{2}\mathcal{W}}{\partial F_{kl}\partial F_{ij}}\left(\textbf{X},\textbf{F}(\textbf{u})\right)\frac{\partial\epsilon_{i}}{\partial X_{j}}\frac{\partial\epsilon_{k}}{\partial X_{l}}dV.
\end{eqnarray*}
Assuming that the following condition holds for all non-zero $\boldsymbol{\epsilon}$ satisfying (\ref{CC:eq:epsilon:1}):
\begin{eqnarray}
\int_{\Omega}\frac{\partial^{2}\mathcal{W}}{\partial F_{kl}\partial F_{ij}}\left(\textbf{X},\textbf{F}(\textbf{u})\right)\frac{\partial\epsilon_{i}}{\partial X_{j}}\frac{\partial\epsilon_{k}}{\partial X_{l}}dV&>&0,\label{CC:eq:ucriterion}
\end{eqnarray}
we obtain (\ref{CC:eq:ubound}), \emph{i.e.} $\textbf{u}\in\mathcal{K}$ is a local minimum for the functional $E_{p}$ defined by (\ref{CC:eq:penergy}). Note that condition (\ref{CC:eq:ucriterion}) is also necessary for the potential energy to attain a local minimum at $\textbf{u}\in\mathcal{K}$.
\eproof

Next, we define the complementary-type strain energy:
\begin{eqnarray*}
\mathcal{W}_{c}(\textbf{F})=\frac{\partial\mathcal{W}}{\partial\textbf{F}}:\textbf{F}-\mathcal{W}(\textbf{F}) &\mbox{in}&\Omega.
\end{eqnarray*}
Then the complementary energy of a stationary admissible field $\textbf{F}'\in\mathcal{S}$ is:
\begin{eqnarray}
E_{c}(\textbf{F}')&=&\int_{\Gamma_{D}}\frac{\partial\mathcal{W}}{\partial\textbf{F}'}\textbf{N}\cdot\textbf{u}_{D}dA-\int_{\Omega}\mathcal{W}_{c}(\textbf{F}')dV.
\label{CC:eq:cenergy}
\end{eqnarray}
In this case also, the contact constraints are present in the set of admissible fields $\mathcal{S}$, but not in the expression of the energy functional $E_{c}$.

The \emph{principle of stationary complementary energy} states that the solution of the CBVP is characterised by the variational formulation:
\[
\mbox{find} \quad \textbf{F}\in \mathcal{S} \quad \mbox{such that}\quad
E_{c}(\textbf{F})=\sup_{\textbf{F}'\in\mathcal{S}}E_{c}(\textbf{F}').
\]

\begin{lemma}
The energy functional $E_{c}$ defined by (\ref{CC:eq:cenergy}) has a local maximum at the solution $\textbf{F}\in\mathcal{S}$ of the complementary-type variational problem (\ref{CC:eq:balance:vgp}), \emph{i.e.}
\begin{eqnarray}
E_{c}(\textbf{F})\geq E_{c}(\textbf{F}'),\label{CC:eq:lbound}
\end{eqnarray}
for all $\textbf{F}'\in\mathcal{S}$.
\end{lemma}

\bproof Let $\textbf{F}'=\textbf{F}+\boldsymbol{\Sigma}\in\mathcal{S}$, where $\textbf{F}\in\mathcal{S}$ is the solution to (\ref{CC:eq:balance:vgp}) and $\Sigma_{ij}\ll 1$, $i,j=1,2,3$, satisfy the following conditions:
\begin{eqnarray}
\frac{\partial}{\partial X_{j}}\left(\frac{\partial^2\mathcal{W}}{\partial F_{ij}\partial F_{kl}}\Sigma_{kl}+\frac{1}{2}\frac{\partial^3\mathcal{W}}{\partial F_{ij}\partial F_{kl}\partial F_{pq}}\Sigma_{pq}\Sigma_{kl}\right)=0 &\mbox{in} & \Omega,\label{CC:eq:sigma:1}\\
\left(\frac{\partial^2\mathcal{W}}{\partial F_{ij}\partial F_{kl}}\Sigma_{kl}+\frac{1}{2}\frac{\partial^3\mathcal{W}}{\partial F_{ij}\partial F_{kl}\partial F_{pq}}\Sigma_{pq}\Sigma_{kl}\right)N_{j}=0 &\mbox{on} & \Gamma_{N}\cup\Gamma_{C}.\label{CC:eq:sigma:2}
\end{eqnarray}
We approximate the value of the complementary energy at $\textbf{F}'$ to second order in $\boldsymbol{\Sigma}$ as follows:
\[
\begin{split}
E_{c}(\textbf{F}')
=&E_{c}(\textbf{F})+\int_{\Gamma_{D}}\left(\frac{\partial\mathcal{W}}{\partial\textbf{F}'}-\frac{\partial\mathcal{W}}{\partial\textbf{F}}\right)\textbf{N}\cdot\textbf{u}_{D}dA
-\int_{\Omega}\left(\frac{\partial\mathcal{W}}{\partial F'_{ij}}-\frac{\partial\mathcal{W}}{\partial F_{ij}}\right)F_{ij}dV\\
&-\int_{\Omega}\left(\frac{\partial\mathcal{W}}{\partial F'_{ij}}-\frac{\partial\mathcal{W}}{\partial F_{ij}}\right)\Sigma_{ij}dV
+\frac{1}{2}\int_{\Omega}\frac{\partial^{2}\mathcal{W}}{\partial F_{ij}\partial F_{kl}}\Sigma_{kl}\Sigma_{ij}dV
\end{split}
\]
By (\ref{CC:eq:balance:vgp}), the above approximation implies:
\begin{eqnarray*}
E_{c}(\textbf{F}')\leq E_{c}(\textbf{F})+\frac{1}{2}\int_{\Omega}\frac{\partial^{2}\mathcal{W}}{\partial F_{ij}\partial F_{kl}}\Sigma_{kl}\Sigma_{ij}dV.
\end{eqnarray*}
Assuming that the following condition is valid for all non-zero $\boldsymbol{\Sigma}$ satisfying (\ref{CC:eq:sigma:1})-(\ref{CC:eq:sigma:2}):
\begin{eqnarray}
\int_{\Omega}\frac{\partial^{2}\mathcal{W}}{\partial F_{ij}\partial F_{kl}}\Sigma_{kl}\Sigma_{ij}dV&=&0,\label{CC:eq:lcriterion}
\end{eqnarray}
we obtain (\ref{CC:eq:lbound}), \emph{i.e.} $\textbf{F}\in\mathcal{S}$ is a local maximum for the functional $E_{c}$ defined by (\ref{CC:eq:cenergy}). The condition (\ref{CC:eq:lcriterion}) is also necessary for the complementary energy to attain a local maximum at $\textbf{F}\in\mathcal{S}$. \eproof

\begin{theorem}
For unconstrained materials, the values of the potential energy functional $E_{p}$ defined by (\ref{CC:eq:penergy}) and of the complementary-type energy functional $E_{c}$ defined by (\ref{CC:eq:cenergy}) represent an upper and a lower bound, respectively, for the minimum potential energy.
\end{theorem}

\bproof
When $\textbf{u}\in\mathcal{K}$ and $\textbf{F}\in\mathcal{S}$ satisfy the CBVP, by (\ref{CC:eq:umixed}):
\begin{eqnarray}
E_{c}(\textbf{F})&=&E_{p}(\textbf{u}).\label{CC:eq:cprinciple}
\end{eqnarray}
Then, by (\ref{CC:eq:ubound}), (\ref{CC:eq:lbound}), and (\ref{CC:eq:cprinciple}):
\begin{equation}\label{CC:eq:cbounds}
E_{c}(\textbf{F}')\leq E_{c}(\textbf{F})=E_{p}(\textbf{u})\leq E_{p}(\textbf{u}'),
\end{equation}
for all $\textbf{u}'\in\mathcal{K}$ and $\textbf{F}'\in\mathcal{S}$.
\eproof

The relation (\ref{CC:eq:cbounds}) guarantees that the values of the potential and complementary energy  provide an upper and a lower bound, respectively, for the (local) minimum potential energy, under the assumption that such a minimum exists.

\section{Materials with Internal Constraints}\label{CC:sec:bounds:con}
In this section, we extend the variational problems formulated in Section~\ref{CC:sec:bounds} to the case where the possible deformations of the body are restricted to those for which:
\begin{equation}\label{CC:eq:constraint}
\gamma(\textbf{F}(\textbf{X}))=0,
\end{equation}
where $\gamma$ is a scalar-valued function representing the material internal constraint and $\textbf{F}(\textbf{X})$ is the deformation gradient at $\textbf{X}$.

\begin{definition} (kinematically admissible fields) The vector field $\textbf{u}$ is a kinematically admissible displacement if and only if it satisfies simultaneously the Dirichlet boundary condition (\ref{CC:eq:Dbc}) and the contact condition (\ref{CC:eq:ndcc}), as well as the internal constraint (\ref{CC:eq:constraint}).
\end{definition}

We denote the closed convex set of kinematically admissible displacement fields by:
\[
\mathcal{K}_{\gamma}=\left\{\textbf{u}'\in W^{1,p}(\Omega;\mathbb{R}^3)\ |\ \gamma\left(\textbf{I}+\mathrm{Grad}\ \textbf{u}'\right)=0\ \mbox{in}\ \Omega,\ \textbf{u}'(\textbf{X})=\textbf{u}_{D} \ \mbox{on} \ \Gamma_{D},\ \eta(\textbf{X}+\textbf{u}'(\textbf{X}))\leq 0\ \mbox{on} \ \Gamma_{C}\right\}.
\]

\begin{definition} (statically admissible fields)
A tensor field $\textbf{F}$ is a statically admissible gradient if and only if it satisfies simultaneously the equilibrium condition (\ref{CC:eq:balance}), the Neumann boundary condition (\ref{CC:eq:Nbc}), the contact conditions (\ref{CC:eq:nscc}) and (\ref{CC:eq:sccc}), and the material constraint (\ref{CC:eq:constraint}).
\end{definition}

We denote the closed convex set of statically admissible gradient fields by:
\[
\begin{split}
\mathcal{S}_{\gamma}=\{(\textbf{F}',\lambda')\in \left[L^{q}(\Omega;\mathbb{R}^3)\right]^2 |\ &\gamma(\textbf{F}')=0,\ \boldsymbol{\mathrm{Div}}\left(\frac{\partial\mathcal{W}}{\partial\textbf{F}'}+\lambda'\frac{\partial\gamma}{\partial\textbf{F}'}\right)=\textbf{0}\ \mbox{in}\ \Omega, \\
&\left(\frac{\partial\mathcal{W}}{\partial\textbf{F}'}+\lambda'\frac{\partial\gamma}{\partial\textbf{F}'}\right)\textbf{N}=\textbf{g}_{N}(\textbf{X}) \ \mbox{on} \ \Gamma_{N}, \\
 &\left(\frac{\partial\mathcal{W}}{\partial\textbf{F}'}+\lambda'\frac{\partial\gamma}{\partial\textbf{F}'}\right)\textbf{N}\cdot\textbf{N}=\textbf{g}_{C}(\textbf{X})\cdot\textbf{N}\leq0\ \mbox{on} \ \Gamma_{C},\\
&\textbf{g}_{C}(\textbf{X})=-\textbf{g}_{C}(\textbf{X}') \ \mbox{if}\ \chi(\textbf{X})=\chi(\textbf{X}') \ \mbox{on} \ \Gamma_{C}\},
\end{split}
\]
where $\lambda'$ is the Lagrange multiplier associated with the material constraint (\ref{CC:eq:constraint}).

Let $\textbf{u}',\textbf{u}''\in\mathcal{K}_{\gamma}$ and $(\textbf{F}'',\lambda'')\in\mathcal{S}_{\gamma}$, such that $\textbf{F}''=\textbf{I}+\mathrm{Grad}\ \textbf{u}''$. The divergence theorem implies:
\begin{equation}\label{CC:eq:mixed:con}
\begin{split}
&\int_{\Gamma_{D}}\left(\frac{\partial\mathcal{W}}{\partial\textbf{F}''}+\lambda''\frac{\partial\gamma}{\partial\textbf{F}''}\right)\textbf{N}\cdot\textbf{u}_{D}dA+\int_{\Gamma_{N}}\textbf{g}_{N}\cdot\textbf{u}'dA+\int_{\Gamma_{C}}\left(\frac{\partial\mathcal{W}}{\partial\textbf{F}''}+\lambda''\frac{\partial\gamma}{\partial\textbf{F}''}\right)\textbf{N}\cdot\textbf{u}'dA\\
&=\int_{\Omega}\mathrm{Div}\left(\frac{\partial\mathcal{W}}{\partial\textbf{F}''}\textbf{u}'+\lambda''\frac{\partial\gamma}{\partial\textbf{F}''}\textbf{u}'\right)dV\\
&=\int_{\Omega}\left(\frac{\partial\mathcal{W}}{\partial\textbf{F}''}+\lambda''\frac{\partial\gamma}{\partial\textbf{F}''}\right):\textbf{F}'dV,
\end{split}
\end{equation}
where $\textbf{F}'=\textbf{I}+\mathrm{Grad}\ \textbf{u}'$. From this identity also, we derive two inequalities.

If, in (\ref{CC:eq:mixed:con}), $\textbf{u}'=\textbf{u}''=\textbf{u}\in\mathcal{K}_{\gamma}$ and $(\textbf{F}'',\lambda'')=(\textbf{F},\lambda)\in\mathcal{S}_{\gamma}$ are the kinematically and statically admissible fields, respectively, satisfying the CBVP with the material constraint (\ref{CC:eq:constraint}), then:
\begin{eqnarray}
\int_{\Gamma_{D}}\left(\frac{\partial\mathcal{W}}{\partial\textbf{F}}+\lambda\frac{\partial\gamma}{\partial\textbf{F}}\right)\textbf{N}\cdot\textbf{u}_{D}dA+\int_{\Gamma_{N}}\textbf{g}_{N}\cdot\textbf{u}dA
&=&\int_{\Omega}\left(\frac{\partial\mathcal{W}}{\partial\textbf{F}}+\lambda\frac{\partial\gamma}{\partial\textbf{F}}\right):\textbf{F}dV.\label{CC:eq:umixed:con}
\end{eqnarray}

\begin{enumerate}
\item First, in (\ref{CC:eq:mixed:con}), let $\textbf{u}''=\textbf{u}\in\mathcal{K}_{\gamma}$ and $(\textbf{F}'',\lambda'')=(\textbf{F},\lambda)\in\mathcal{S}_{\gamma}$ satisfy the CBVP with the material constraint (\ref{CC:eq:constraint}), and $\textbf{u}'\in\mathcal{K}_{\gamma}$. Subtracting (\ref{CC:eq:umixed:con}) from (\ref{CC:eq:mixed:con}) gives:
\begin{eqnarray*}
\int_{\Gamma_{N}}\textbf{g}_{N}\cdot(\textbf{u}'-\textbf{u})dA+\int_{\Gamma_{C}}\left(\frac{\partial\mathcal{W}}{\partial\textbf{F}}+\lambda\frac{\partial\gamma}{\partial\textbf{F}}\right)\textbf{N}\cdot\textbf{u}'dA
&=&\int_{\Omega}\left(\frac{\partial\mathcal{W}}{\partial\textbf{F}}+\lambda\frac{\partial\gamma}{\partial\textbf{F}}\right):(\textbf{F}'-\textbf{F})dV.
\end{eqnarray*}

In this case, the primal variational problem is to find $\textbf{u}\in\mathcal{K}_{\gamma}$ that satisfies:
\begin{equation}\label{CC:eq:balance:vdp:con}
\int_{\Omega}\left[\frac{\partial\mathcal{W}}{\partial\textbf{F}}\left(\textbf{X},\textbf{F}(\textbf{u})\right)+\lambda\frac{\partial\gamma}{\partial\textbf{F}}\left(\textbf{F}(\textbf{u})\right)\right]: \mathrm{Grad}(\textbf{u}'-\textbf{u})dV
\geq\int_{\Gamma_{N}}\textbf{g}_{N}\cdot(\textbf{u}'-\textbf{u})dA,
\end{equation}
for all kinematically admissible fields $\textbf{u}'\in\mathcal{K}_{\gamma}$.

\item Next, in (\ref{CC:eq:mixed:con}), let $\textbf{u}'=\textbf{u}\in\mathcal{K}_{\gamma}$ and $(\textbf{F}',\lambda')=(\textbf{F},\lambda)\in\mathcal{S}_{\gamma}$ satisfy the CBVP with the material constraint (\ref{CC:eq:constraint}), and $(\textbf{F}'',\lambda'')\in\mathcal{S}_{\gamma}$. Subtracting (\ref{CC:eq:umixed:con}) from (\ref{CC:eq:mixed:con}) yields:
\[
\begin{split}
&\int_{\Gamma_{D}}\left(\frac{\partial\mathcal{W}}{\partial\textbf{F}''}+\lambda''\frac{\partial\gamma}{\partial\textbf{F}''}-\frac{\partial\mathcal{W}}{\partial\textbf{F}}-\lambda\frac{\partial\gamma}{\partial\textbf{F}}\right)\textbf{N}\cdot\textbf{u}_{D}dA\\
&+\int_{\Gamma_{C}}\left(\frac{\partial\mathcal{W}}{\partial\textbf{F}''}+\lambda''\frac{\partial\gamma}{\partial\textbf{F}''}\right)\textbf{N}\cdot\textbf{u}dA\\
&=\int_{\Omega}\left(\frac{\partial\mathcal{W}}{\partial\textbf{F}''}+\lambda''\frac{\partial\gamma}{\partial\textbf{F}''}-\frac{\partial\mathcal{W}}{\partial\textbf{F}}-\lambda\frac{\partial\gamma}{\partial\textbf{F}}\right):\textbf{F}dV.
\end{split}
\]

Then the complementary-type variational problem is to find $(\textbf{F},\lambda)\in\mathcal{S}_{\gamma}$ satisfying:
\begin{equation}\label{CC:eq:balance:vgp:con}
\int_{\Omega}\left(\frac{\partial\mathcal{W}}{\partial\textbf{F}'}+\lambda'\frac{\partial\gamma}{\partial\textbf{F}'}-\frac{\partial\mathcal{W}}{\partial\textbf{F}}-\lambda\frac{\partial\gamma}{\partial\textbf{F}}\right):\textbf{F}dV
\geq\int_{\Gamma_{D}}\left(\frac{\partial\mathcal{W}}{\partial\textbf{F}'}+\lambda'\frac{\partial\gamma}{\partial\textbf{F}'}-\frac{\partial\mathcal{W}}{\partial\textbf{F}}-\lambda\frac{\partial\gamma}{\partial\textbf{F}}\right)\textbf{N}\cdot\textbf{u}_{D}dA,
\end{equation}
for all statically admissible fields $(\textbf{F}',\lambda')\in\mathcal{S}_{\gamma}$.
\end{enumerate}

\subsection{Upper and Lower Bounds on the Strain Energy}
In this case, the potential energy of a kinematically admissible field $\textbf{u}'\in\mathcal{K}_{\gamma}$ is given by:
\begin{eqnarray}
E_{p}(\textbf{u}')&=&\int_{\Omega}\mathcal{W}\left(\textbf{X},\textbf{F}'(\textbf{u}')\right)dV-\int_{\Gamma_{N}}\textbf{g}_{N}\cdot\textbf{u}'dA.\label{CC:eq:penergy:con}
\end{eqnarray}

\begin{lemma}\label{CC:lemma:Ep:con}
The energy functional $E_{p}$ defined by (\ref{CC:eq:penergy:con}) has a local minimum at the solution $\textbf{u}\in\mathcal{K}_{\gamma}$ of (\ref{CC:eq:balance:vdp:con}), \emph{i.e.}
\begin{eqnarray}
E_{p}(\textbf{u})&\leq&E_{p}(\textbf{u}'),\label{CC:eq:ubound:con}
\end{eqnarray}
for all $\textbf{u}'\in\mathcal{K}_{\gamma}$.
\end{lemma}

\bproof Let $\textbf{u}'=\textbf{u}+\boldsymbol{\epsilon}\in\mathcal{K}_{\gamma}$, where $\textbf{u}\in\mathcal{K}_{\gamma}$ is the solution to (\ref{CC:eq:balance:vdp:con}) and $\epsilon_{i}\ll 1$, $i=1,2,3$, such that:
\begin{eqnarray}
\frac{\partial\gamma}{\partial F_{ij}}\frac{\partial\epsilon_{i}}{\partial X_{j}}=0\ \mbox{in}\ \Omega,\label{CC:eq:epsilon:0:con}\\
\boldsymbol{\epsilon}=0\ \mbox{on}\  \Gamma_{D}\cup\Gamma_{C}.\label{CC:eq:epsilon:1:con}
\end{eqnarray}
We approximate the value of the potential energy at $\textbf{u}'$ to the second order in $\boldsymbol{\epsilon}$ as follows:
\[
\begin{split}
E_{p}(\textbf{u}',\lambda')
&=E_{p}(\textbf{u})
+\int_{\Omega}\left[\frac{\partial\mathcal{W}}{\partial F_{pq}}\left(\textbf{X},\textbf{F}(\textbf{u})\right)+\lambda\frac{\partial\gamma}{\partial F_{pq}}\left(\textbf{F}(\textbf{u})\right)\right]\frac{\partial\epsilon_{p}}{\partial X_{q}}dV-\int_{\Gamma_{N}}\textbf{g}_{N}\cdot\boldsymbol{\epsilon}dA\\
&+\frac{1}{2}\int_{\Omega}\left[\frac{\partial^{2}\mathcal{W}}{\partial F_{kl}\partial F_{pq}}\left(\textbf{X},\textbf{F}(\textbf{u})\right)+\lambda\frac{\partial^{2}\gamma}{\partial F_{kl}\partial F_{pq}}\left(\textbf{F}(\textbf{u})\right)\right]\frac{\partial\epsilon_{p}}{\partial X_{q}}\frac{\partial\epsilon_{k}}{\partial X_{l}}dV.
\end{split}
\]
By (\ref{CC:eq:balance:vdp:con}):
\begin{eqnarray*}
\int_{\Omega}\left[\frac{\partial\mathcal{W}}{\partial F_{pq}}\left(\textbf{X},\textbf{F}(\textbf{u})\right)+\lambda\frac{\partial\gamma}{\partial F_{pq}}\left(\textbf{F}(\textbf{u})\right)\right]\frac{\partial\epsilon_{p}}{\partial X_{q}}dV&\geq&\int_{\Gamma_{N}}\textbf{g}_{N}\cdot\boldsymbol{\epsilon}dA.
\end{eqnarray*}
Hence:
\begin{eqnarray*}
E_{p}(\textbf{u}')&\geq&E_{p}(\textbf{u})+
\frac{1}{2}\int_{\Omega}\left[\frac{\partial^{2}\mathcal{W}}{\partial F_{kl}\partial F_{pq}}\left(\textbf{X},\textbf{F}(\textbf{u})\right)+\lambda\frac{\partial^{2}\gamma}{\partial F_{kl}\partial F_{pq}}\left(\textbf{F}(\textbf{u})\right)\right]\frac{\partial\epsilon_{p}}{\partial X_{q}}\frac{\partial\epsilon_{k}}{\partial X_{l}}dV.
\end{eqnarray*}
Assuming that the following condition holds for all non-zero $\boldsymbol{\epsilon}$ satisfying (\ref{CC:eq:epsilon:0:con})- (\ref{CC:eq:epsilon:1:con}):
\begin{eqnarray}
\int_{\Omega}\left[\frac{\partial^{2}\mathcal{W}}{\partial F_{kl}\partial F_{pq}}\left(\textbf{X},\textbf{F}(\textbf{u})\right)+\lambda\frac{\partial^{2}\gamma}{\partial F_{kl}\partial F_{pq}}\left(\textbf{F}(\textbf{u})\right)\right]\frac{\partial\epsilon_{p}}{\partial X_{q}}\frac{\partial\epsilon_{k}}{\partial X_{l}}dV&>&0,\label{CC:eq:ucriterion:con}
\end{eqnarray}
we obtain (\ref{CC:eq:ubound:con}), \emph{i.e.} $\textbf{u}\in\mathcal{K}_{\gamma}$ is a local minimum for the functional $E_{p}$ defined by (\ref{CC:eq:penergy:con}). \eproof

We also define the complementary-type energy density as follows:
\begin{eqnarray}
\mathcal{W}_{c}(\textbf{F},\lambda)&=&\left(\frac{\partial\mathcal{W}}{\partial\textbf{F}}+\lambda\frac{\partial\gamma}{\partial\textbf{F}}\right):\textbf{F}-\mathcal{W}(\textbf{F}).
\end{eqnarray}
The corresponding complementary energy of a statically admissible field $(\textbf{F}',\lambda')\in\mathcal{S}_{\gamma}$ is:
\begin{eqnarray}
E_{c}(\textbf{F}',\lambda')&=&\int_{\Gamma_{D}}\left(\frac{\partial\mathcal{W}}{\partial\textbf{F}'}+\lambda'\frac{\partial\gamma}{\partial\textbf{F}'}\right)\textbf{N}\cdot\textbf{u}_{D}dA-\int_{\Omega}\mathcal{W}_{c}(\textbf{F}',\lambda')dV.\label{CC:eq:cenergy:con}
\end{eqnarray}

\begin{lemma}\label{CC:lemma:Ec:con}
The energy functional $E_{c}$ defined by (\ref{CC:eq:cenergy:con}) has a local maximum at the solution $(\textbf{F},\lambda)\in\mathcal{S}_{\gamma}$ of (\ref{CC:eq:balance:vgp:con}), \emph{i.e.}
\begin{eqnarray}
E_{c}(\textbf{F},\lambda)&\geq&E_{c}(\textbf{F}',\lambda'),\label{CC:eq:lbound:con}
\end{eqnarray}
for all $(\textbf{F}',\lambda')\in\mathcal{S}_{\gamma}$.
\end{lemma}

\bproof Let $(\textbf{F}',\lambda')\in\mathcal{S}_{\gamma}$, such that $\textbf{F}'=\textbf{F}+\boldsymbol{\Sigma}$ and $\lambda'=\lambda+\epsilon$, where $(\textbf{F},\lambda)\in\mathcal{S}_{\gamma}$ is the solution to (\ref{CC:eq:balance:vgp:con}), and $\Sigma_{ij}\ll 1$, $i,j=1,2,3$, and $\epsilon\ll 1$ satisfy the following conditions:
\begin{eqnarray}
\frac{\partial\gamma}{\partial F_{ij}}\Sigma_{ij}=0 & \mbox{in} & \Omega,\label{CC:eq:sigma:0:con}\\
\frac{\partial}{\partial X_{j}}\left(\frac{\partial^2\mathcal{W}}{\partial F_{ij}\partial F_{kl}}\Sigma_{kl}+\epsilon\frac{\partial\gamma}{\partial F_{ij}}+\lambda\frac{\partial^2\gamma}{\partial F_{ij}\partial F_{kl}}\Sigma_{kl}\right)=0 & \mbox{in} & \Omega,\label{CC:eq:sigma:1:con}\\
\left(\frac{\partial^2\mathcal{W}}{\partial F_{ij}\partial F_{kl}}\Sigma_{kl}+\epsilon\frac{\partial\gamma}{\partial F_{ij}}+\lambda\frac{\partial^2\gamma}{\partial F_{ij}\partial F_{kl}}\Sigma_{kl}\right)N_{j}=0 & \mbox{on} & \Gamma_{C}\cup\Gamma_{N}.\label{CC:eq:sigma:2:con}
\end{eqnarray}
We approximate the value of the complementary energy at $\textbf{F}'$ to second order in $\boldsymbol{\Sigma}$ and $\epsilon$ as follows:
\[
\begin{split}
E_{c}(\textbf{F}',\lambda')
&=E_{c}(\textbf{F},\lambda)+\int_{\Gamma_{D}}\left(\frac{\partial\mathcal{W}}{\partial\textbf{F}'}+\lambda'\frac{\partial\gamma}{\partial\textbf{F}'}-\frac{\partial\mathcal{W}}{\partial\textbf{F}}-\lambda\frac{\partial\gamma}{\partial\textbf{F}}\right)\textbf{N}\cdot\textbf{u}_{D}dA\\
&-\int_{\Omega}\left(\frac{\partial\mathcal{W}}{\partial F'_{ij}}+\lambda'\frac{\partial\gamma}{\partial F'_{ij}}-\frac{\partial\mathcal{W}}{\partial F_{ij}}-\lambda\frac{\partial\gamma}{\partial F_{ij}}\right)F_{ij}dV\\
&-\int_{\Omega}\left(\frac{\partial\mathcal{W}}{\partial F'_{ij}}+\lambda'\frac{\partial\gamma}{\partial F'_{ij}}-\frac{\partial\mathcal{W}}{\partial F_{ij}}-\lambda\frac{\partial\gamma}{\partial F_{ij}}\right)\Sigma_{ij}dV\\
&+\frac{1}{2}\int_{\Omega}\left(\frac{\partial^{2}\mathcal{W}}{\partial F_{ij}\partial F_{kl}}+\lambda\frac{\partial^2\gamma}{\partial F_{ij}\partial F_{kl}}\right)\Sigma_{kl}\Sigma_{ij}dV
\end{split}
\]
By (\ref{CC:eq:balance:vgp:con}), the above approximation implies:
\begin{eqnarray*}
E_{c}(\textbf{F}',\lambda')\leq E_{c}(\textbf{F},\lambda)+\frac{1}{2}\int_{\Omega}\left(\frac{\partial^{2}\mathcal{W}}{\partial F_{ij}\partial F_{kl}}+\lambda\frac{\partial^2\gamma}{\partial F_{ij}\partial F_{kl}}\right)\Sigma_{kl}\Sigma_{ij}dV.
\end{eqnarray*}
Assuming that the following condition is valid for all non-zero $\boldsymbol{\Sigma}$ and $\boldsymbol{\epsilon}$ satisfying (\ref{CC:eq:sigma:0:con})-(\ref{CC:eq:sigma:2:con}):
\begin{eqnarray}
\int_{\Omega}\left(\frac{\partial^{2}\mathcal{W}}{\partial F_{ij}\partial F_{kl}}+\lambda\frac{\partial^2\gamma}{\partial F_{ij}\partial F_{kl}}\right)\Sigma_{kl}\Sigma_{ij}dV&>&0,\label{CC:eq:lcriterion:con}
\end{eqnarray}
we obtain (\ref{CC:eq:lbound:con}), \emph{i.e.} $(\textbf{F},\lambda)\in\mathcal{S}_{\gamma}$ represents a local maximum for the functional $E_{c}$ defined by (\ref{CC:eq:cenergy:con}). \eproof

\begin{theorem}
For materials with the internal constraint (\ref{CC:eq:constraint}), the values of the potential energy $E_{p}$ defined by (\ref{CC:eq:penergy:con}) and of the complementary type energy $E_{c}$ defined (\ref{CC:eq:cenergy:con}) provide an upper and a lower bound, respectively, for the minimum potential energy.
\end{theorem}

\bproof
When $\textbf{u}\in\mathcal{K}_{\gamma}$ and $(\textbf{F},\lambda)\in\mathcal{S}_{\gamma}$ satisfy the CBVP with the material constraint (\ref{CC:eq:constraint}), by (\ref{CC:eq:umixed:con}):
\begin{eqnarray}
E_{c}(\textbf{F},\lambda)&=&E_{p}(\textbf{u}).\label{CC:eq:cprinciple:con}
\end{eqnarray}
Then (\ref{CC:eq:ubound:con}), (\ref{CC:eq:lbound:con}), and (\ref{CC:eq:cprinciple:con}) imply:
\begin{equation}\label{CC:eq:cbounds:con}
E_{c}(\textbf{F}',\lambda')\leq E_{c}(\textbf{F},\lambda)=E_{p}(\textbf{u})\leq E_{p}(\textbf{u}'),
\end{equation}
for all $\textbf{u}'\in\mathcal{K}_{\gamma}$ and $(\textbf{F}',\lambda')\in\mathcal{S}_{\gamma}$. \eproof\\

The relation (\ref{CC:eq:cbounds:con}) guarantees that the potential and complementary energy  form an enclosure for the (local) minimum potential energy, under the assumption that such a minimum exists.

\section{Extension to Two Elastic Bodies in Unilateral Contact}\label{CC:sec:problem:two}

\begin{figure}[htbp]
\begin{center}
\scalebox{0.5}{\includegraphics{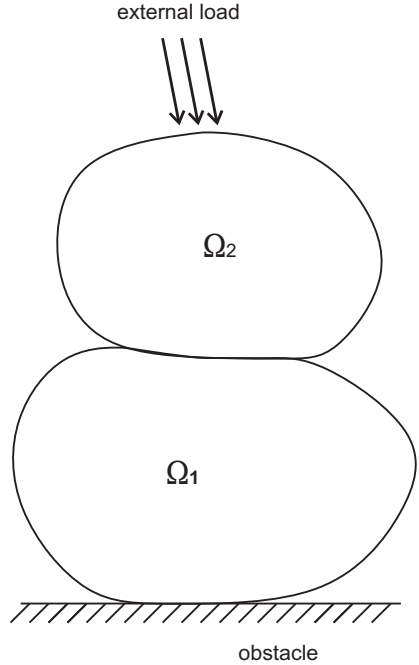}}
\caption{Schematic representation of a system of two elastic bodies in mutual non-penetrative contact and subject to external loading.}
\label{fig:contact-two}
\end{center}
\end{figure}

We now consider a system formed from two elastic bodies made from possibly different homogeneous hyperelastic materials described by the strain energy function $\mathcal{W}_{i}$, $i=1,2$, respectively, which are in mutual non-penetrative contact on part of their boundary, and denote by  $\Omega_{i}$, $i=1,2$, the two open, bounded, connected, with Lipschitz continuous boundary, distinct domains occupied by the two bodies, respectively, such that $\Omega_{1}\cap\Omega_{2}=\emptyset$ (see Figure~\ref{fig:contact-two}).

Each body is subject to a finite elastic deformation:
\[
\boldsymbol{\chi}_{i}:\Omega_{i}\to\mathbb{R}^{3},
\]
such that $\boldsymbol{\chi}_{i}$ satisfies (\ref{CC:eq:inj1})-(\ref{CC:eq:inj2}) on $\Omega_{i}$.

The corresponding Lagrangian equation is:
\begin{eqnarray*}
\boldsymbol{\mathrm{Div}}\textbf{P}(\textbf{X})=\textbf{0} &\mbox{in}& \Omega_{1}\cup\Omega_{2},
\end{eqnarray*}
such that
\begin{eqnarray*}
\textbf{P}=\frac{\partial\mathcal{W}_{i}}{\partial\textbf{F}_{i}}-p\textbf{F}_{i}^{-T} &\mbox{in}& \Omega_{i},\ i=1,2,
\end{eqnarray*}
where $\textbf{F}_{i}=\mathrm{Grad}\ \boldsymbol{\chi}_{i}$, $i=1,2$, $p=0$ for compressible materials, and $p$ is the hydrostatic pressure for incompressible materials.

For this system, boundary conditions similar to (\ref{CC:eq:Dbc})-(\ref{CC:eq:Nbc}) can be imposed on the external surface of individual bodies, while the conditions of non-penetrative contact on the surface $\Gamma_{C}$ take the form (\ref{CC:eq:ndcc})-(\ref{CC:eq:sccc}), where, for mutual or self-contact, $\eta=[\textbf{u}(\textbf{X})]\cdot\textbf{N}$ where $[\textbf{u}(\textbf{X})]=\textbf{u}(\textbf{X})-\textbf{u}(\textbf{X}')$ denotes the jump in the displacement at the points $\textbf{X}\neq\textbf{X}'$ where contact may occur and $\textbf{N}$ and $\textbf{N}'$ are the outward unit normal vectors to the boundary at those points, respectively.

The formulation of the corresponding variational problems and their analysis are then parallel to that presented in Sections~\ref{CC:sec:bounds}-\ref{CC:sec:bounds:con}. We illustrate this with the following examples for large elastic deformations which can be maintained in every homogeneous isotropic incompressible material in the absence of body forces.

\begin{example}
Two thin elastic bodies which occupy the domains $\Omega_{1}=(0,1/2)\times (0,d)\times(0,d)$ and $\Omega_{2}=(1/2,1)\times (0,d)\times(0,d)$, respectively, where $1\ll d<+\infty$, are made from different neo-Hookean materials described by the strain energy function:
\begin{equation}\label{CC:eq:NH}
\mathcal{W}_{i}=\frac{C_{i}}{2}(I_{1}-3), \qquad i=1,2,
\end{equation}
where $C_{i}$, $i=1,2$, are the corresponding material constants and $I_{1}$ is the first principal strain invariant. Each body occupying the domain $\Omega_{i}$, $i=1,2$, is subject to a triaxial stretch of the form:
\begin{equation}\label{CC:eq:stretch}
x=a_{i}X+b_{i},\qquad y=\frac{Y}{\sqrt{a_{i}}},\qquad z=\frac{Z}{\sqrt{a_{i}}}, \qquad i=1,2,
\end{equation}
where $a_{i}$, $b_{i}$, $i=1,2$, are positive constants, and $(X,Y,Z)$ and $(x,y,z)$ are the Cartesian coordinates for the reference and the deformed configuration, respectively. During the deformation, non-penetrative contact between the two bodies at the interface $\Gamma_{C}=\{1/2\}\times (0,d)\times(0,d)$ is assumed. Since the dimension of the bodies in the $X$-direction is much smaller than in the other two directions, changes in the area of the surface $\Gamma_{C}$ during deformation can be neglected (this is indicated by the dashed line at the horizontal ends of the domains in Figure~\ref{fig:compress-two}). This deformation is attained by uniformly loading the surface $\Gamma_{N}=\{0\}\times (0,d)\times(0,d)$ of the first body, while on the surface $\Gamma_{D}=\{1\}\times (0,d)\times(0,d)$ of the second body the deformation is prescribed.

\begin{figure}[htbp]
\begin{center}
\scalebox{0.5}{\includegraphics{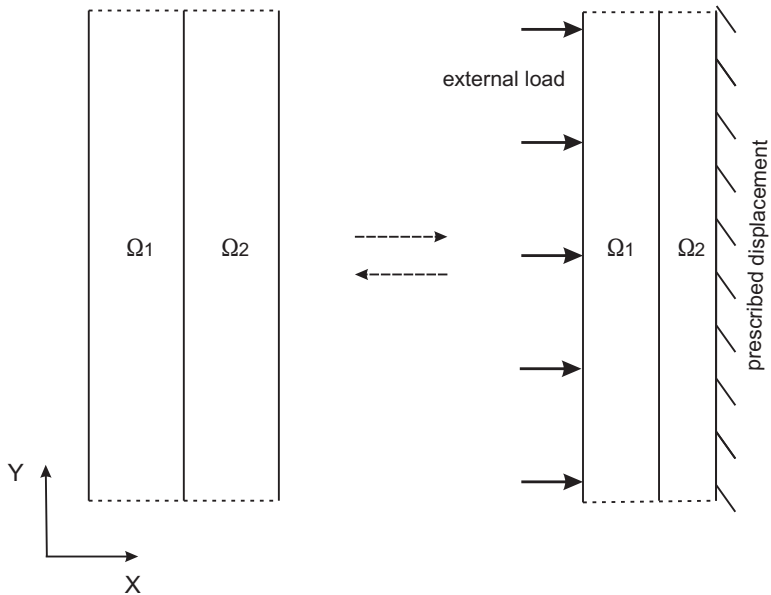}}
\caption{Schematic representation of the undeformed (left) and deformed (right) states of a system of two thin elastic bodies in mutual non-penetrative contact, compressing under external load.}
\label{fig:compress-two}
\end{center}
\end{figure}

For the deformation (\ref{CC:eq:stretch}), the $x$, $y$ and $z$-directions are principal directions, and the principal stretches are respectively:
\[
\lambda_{1}=a_{i},\qquad \lambda_{2}=\frac{1}{\sqrt{a_{i}}},\qquad \lambda_{3}=\frac{1}{\sqrt{a_{i}}}.
\]
We first account for the unilateral contact constraints at the interface $\Gamma_{C}$, where the relative normal displacement must satisfy the non-penetration condition (\ref{CC:eq:ndcc}):
\begin{equation}\label{CC:eq:nonpenetration:1}
\frac{a_{1}}{2}+b_{1}\leq \frac{a_{2}}{2}+b_{2}.
\end{equation}
Since, by the action-reaction law the contact stresses at $\Gamma_{C}$ are equal in magnitude, if the hydrostatic pressure for the body occupying the domain $\Omega_{i}$ is equal to $-p_{i}$, $i=1,2$, then by the normal contact stresses condition (\ref{CC:eq:nscc}):
\begin{equation}\label{CC:eq:compression:1}
-p_{1}+C_{1}a_{1}^2=-p_{2}+C_{2}a_{2}^2\leq 0,
\end{equation}
and by the complementarity condition (\ref{CC:eq:cpcc}):
\begin{equation}\label{CC:eq:complementarity:1}
\left(\frac{a_{1}}{2}+b_{1}-\frac{a_{2}}{2}-b_{2}\right)\left(-p_{1}+C_{1}a_{1}^2\right)= 0.
\end{equation}
When the equality holds in (\ref{CC:eq:nonpenetration:1}), the conditions (\ref{CC:eq:ucriterion:con}) and (\ref{CC:eq:lcriterion:con}) are satisfied simultaneously if and only if:
\begin{equation}\label{CC:eq:bounds:1}
-C_{i}\sqrt{a_{i}}<p_{i}<C_{i}\sqrt{a_{i}},\qquad i=1,2,
\end{equation}
such that (\ref{CC:eq:compression:1}) is valid. Equivalently, the external load $\tau$ in the $x$-direction per unit area of the deformed configuration satisfies:
\[
-\min\left\{C_{1}\left(\sqrt{a_{1}}-a_{1}^2\right),C_{2}\left(\sqrt{a_{2}}-a_{2}^2\right)\right\}<\tau< 0,
\]
and the action-reaction equality in (\ref{CC:eq:compression:1}) holds. If the strict inequality is satisfied in (\ref{CC:eq:nonpenetration:1}), then $\tau=0$.
\end{example}

\begin{example}
For the two elastic bodies described in the previous example and deformed by (\ref{CC:eq:stretch}), we now consider the case where cohesive contact occurs at $\Gamma_{C}$. Then the relation (\ref{CC:eq:compression:1}) is replaced by:
\begin{equation}\label{CC:eq:cohesion:1}
-p_{1}+C_{1}a_{1}^2=-p_{2}+C_{2}a_{2}^2\leq g,
\end{equation}
where $g>0$ is the given cohesion parameter, and the complementarity condition (\ref{CC:eq:complementarity:1}) becomes:
\begin{equation}\label{CC:eq:complementarity:g1}
\left(\frac{a_{1}}{2}+b_{1}-\frac{a_{2}}{2}-b_{2}\right)\left(-p_{1}+C_{1}a_{1}^2-g\right)= 0.
\end{equation}
In this case, if the equality is valid in (\ref{CC:eq:nonpenetration:1}), then by (\ref{CC:eq:bounds:1}) and (\ref{CC:eq:cohesion:1}), the external load $\tau$ in the $x$-direction per unit area of the deformed configuration satisfies:
\[
-\min_{i=1,2}C_{i}\left(\sqrt{a_{i}}-a_{i}^2\right)<\tau
<\min\left\{g,\min_{i=1,2}C_{i}\left(\sqrt{a_{i}}+a_{i}^2\right)\right\},
\]
and the action-reaction equality in (\ref{CC:eq:cohesion:1}) holds. If the strict inequality is valid in (\ref{CC:eq:nonpenetration:1}), then $\tau=g$.
\end{example}

\begin{example}
We consider again two thin bodies made from different neo-Hookean materials as described by (\ref{CC:eq:NH}) and occupying the domains $\Omega_{1}=(0,1/2)\times (0,d)\times(0,d)$ and $\Omega_{2}=(1/2,1)\times (0,d)\times(0,d)$, respectively, where $1\ll d<+\infty$. Each body occupying the domain $\Omega_{i}$, $i=1,2$, is now subject to the a combined stretch and bending of the form:
\begin{equation}\label{CC:eq:bend}
r=\sqrt{2a_{i}X+b_{i}},\qquad \theta=\frac{AY}{\sqrt{a}_{i}},\qquad z=\frac{Z}{A\sqrt{a}_{i}},
\end{equation}
where $A$ and $a_{i}$, $b_{i}$, $i=1,2$ are positive constants, $(X,Y,Z)$ are the Cartesian coordinates for the reference configuration and $(r,\theta,z)$ are the cylindrical polar coordinates for the deformed configuration, such that $r\in(r_{i-1},r_{i})$, $i=1,2$. During the deformation, non-penetrative contact between the two bodies at the interface $\Gamma_{C}=\{1/2\}\times (0,d)\times(0,d)$ is assumed (see Figure~\ref{fig:bend-two}). This deformation is attained by uniformly loading the surface $\Gamma_{N}=\{0\}\times (0,d)\times(0,d)$ of the first body, while on the surface $\Gamma_{D}=\{1\}\times (0,d)\times(0,d)$ of the second body the deformation is prescribed. In this case also, since the dimension of the bodies in the $X$-direction is much smaller than in the other two directions, changes in the area of the contact surface during deformation may be neglected.

\begin{figure}[htbp]
\begin{center}
\scalebox{0.5}{\includegraphics{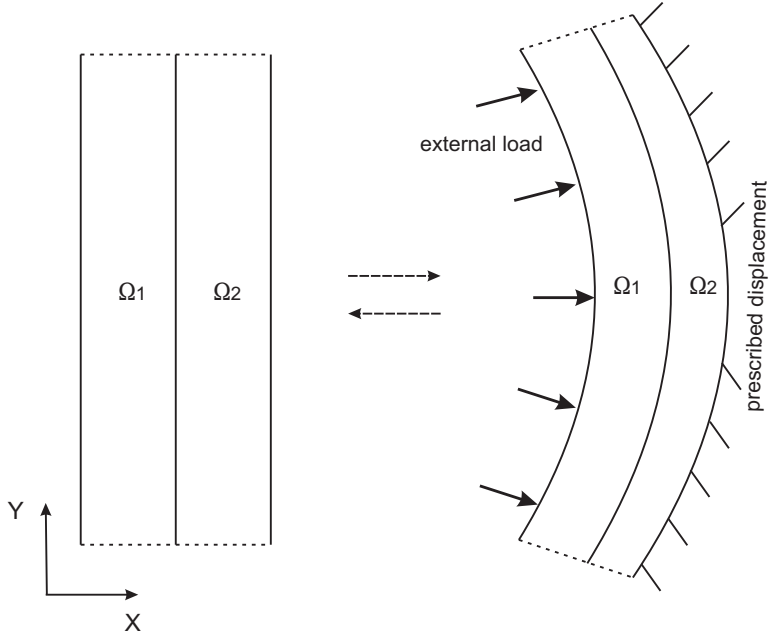}}
\caption{Schematic representation of the undeformed (left) and deformed (right) states of a system of two thin elastic bodies in mutual non-penetrative contact, bending under external load.}
\label{fig:bend-two}
\end{center}
\end{figure}

For the deformation  (\ref{CC:eq:bend}), the $r$, $\theta$ and $z$-directions are principal directions, and the principal stretches are:
\[
\lambda_{1}=\frac{a_{i}}{r},\qquad \lambda_{2}=\frac{Ar}{\sqrt{a_{i}}},\qquad \lambda_{3}=\frac{1}{A\sqrt{a_{i}}}.
\]
At the contact interface $\Gamma_{C}$, the relative normal displacement satisfies the non-penetration condition:
\begin{equation}\label{CC:eq:nonpenetration:2}
\sqrt{a_{1}+b_{1}}\leq\sqrt{a_{2}+b_{2}},
\end{equation}
and setting the hydrostatic pressure $-p_{i}$ for the body occupying the domain $\Omega_{i}$, $i=1,2$, the normal contact stresses satisfy:
\begin{equation}\label{CC:eq:compression:2}
-p_{1}+C_{1}\frac{a_{1}^2}{r_{1}^2}=-p_{2}+C_{2}\frac{a_{2}^2}{r_{1}^2}\leq 0.
\end{equation}
By the complementarity condition (\ref{CC:eq:cpcc}):
\begin{equation}\label{CC:eq:complementarity:2}
\left(a_{1}+b_{1}-a_{2}-b_{2}\right)\left(-p_{1}+C_{1}\frac{a_{1}^2}{r_{1}^2}\right)= 0.
\end{equation}
When the equality holds in (\ref{CC:eq:nonpenetration:1}), the conditions (\ref{CC:eq:ucriterion:con}) and (\ref{CC:eq:lcriterion:con}) are satisfied simultaneously if and only if:
\[
-\frac{C_{i}}{\max\{a_{i}/r_{i-1},Ar_{i}/\sqrt{a_{i}},1/(A\sqrt{a_{i}})\}}<p_{i}<\frac{C_{i}}{\max\{a_{i}/r_{i-1},Ar_{i}/\sqrt{a_{i}},1/(A\sqrt{a_{i}})\}}, \qquad i=1,2,
\]
such that (\ref{CC:eq:compression:2}) is valid. Equivalently, the external load $\tau$ in the $r$-direction per unit area of the deformed configuration satisfies:
\[
-\min\left\{\frac{C_{1}}{\max\{a_{1}/r_{0},Ar_{1}/\sqrt{a_{1}},1/(A\sqrt{a_{1}})\}}-\frac{C_{1}a_{1}^2}{r_{1}^2},\frac{C_{2}}{\max\{a_{2}/r_{1},Ar_{2}/\sqrt{a_{2}},1/(A\sqrt{a_{2}})\}}-\frac{C_{2}a_{2}^2}{r_{1}^2}\right\}<\tau<0,
\]
and the action-reaction equality in (\ref{CC:eq:compression:2}) holds. If the strict inequality is satisfied in (\ref{CC:eq:nonpenetration:2}), then $\tau=0$.
\end{example}

\section{Conclusion}\label{CC:sec:conclusion}
In large deformations, elastic bodies are very likely to enter into contact with neighbouring obstacles or with themselves, and therefore both the mathematical models and their analysis must take these phenomena into account. In this paper, for an elastic body made from an homogeneous isotropic non-linear hyperleastic material in unilateral contact with a rigid obstacle or with another elastic body, trial functions for the deformation gradient were used to obtain a variational principle of stationary complementary energy type. Then the stationary potential energy and complementary energy principles were used to provide an enclosure on the total strain energy of the finitely deformed body. To illustrate the theory, this variational framework was applied to determine upper and lower bounds on the external load for a system made from two elastic bodies of neo-Hookean material under uniform compression or under combined stretch and bending. This investigation addresses the need for a better understanding of contact problems in natural and industrial systems where mechanical models that take into account the large stresses and strains at adjoining material surfaces are crucial.

\section*{Acknowledgements}
The support for L.A.M. by the Engineering and Physical Sciences Research Council of Great Britain under research grant EP/M011992/1 is gratefully acknowledged. There are no data associated with this paper.


\end{document}